\documentclass[10pt]{amsart}

\theoremstyle{plain}
\newtheorem{Theorem}{Theorem}

\newtheorem{Lemma}[Theorem]{Lemma}
\newtheorem{Proposition}[Theorem]{Proposition}

\theoremstyle{definition}

\theoremstyle{remark}

\DeclareMathOperator{\codim}{codim}
\DeclareMathOperator{\discrep}{discrep}

\DeclareMathOperator{\divisor}{div}
\DeclareMathOperator{\Exc}{Exc}
\DeclareMathOperator{\md}{md}

\DeclareMathOperator{\ord}{ord}

\DeclareMathOperator{\Reg}{Reg}
\DeclareMathOperator{\Supp}{Supp}

\newcommand{\Q}{\ensuremath{\mathbb{Q}\,}}

\newcommand{\C}{\ensuremath{\mathbb{C}\,}}
\newcommand{\A}{\ensuremath{\mathbb{A}\,}}

\newcommand{\Ocal}{\mathcal{O}}
\newcommand{\Mcal}{\mathcal{M}}

\newcommand{\bbu}{\mathbf{u}}
\newcommand{\bby}{\mathbf{y}}

\newcommand{\union}{\cup}
\newcommand{\inters}{\cap}

\newcommand{\notsubset}{\not\subseteq}

\begin{document}

%%%%%%%%%%%%%%%%%%%%%%%%%%%%%%%%%%%%%%%%%%%%%%%%%%%%%%%%%%%%
%%%%%%%%%%%%%%%%%%%%%%%%%%%%%%%%%%%%%%%%%%%%%%%%%%%%%%%%%%%%

%%% BEGIN TOPMATTER

    % __TITLE__
\title[Minimal discrepancies of hypersurfaces]
		{Minimal discrepancies of hypersurface singularities}

    % __AUTHOR__
\author{Vladimir Ma\c{s}ek}
\address{Department of Mathematics, Box 1146, Washington University,
		St. Louis, MO 63130}
\email{vmasek@math.wustl.edu}

    % __MATH. SUBJ. CLASSIF. (AMS, 1991)__
\subjclass{Primary 14E15; Secondary 14B05, 14E30, 14J70}

    % __ABSTRACT__
\begin{abstract}

We give an upper bound for the minimal discrepancies of hypersurface
singularities. As an application, we show that Shokurov's conjecture is true
for log-terminal threefolds.

\end{abstract}

\maketitle

%%% END TOPMATTER

%%%%%%%%%%%%%%%%%%%%%%%%%%%%%%%%%%%%%%%%%%%%%%%%%%%%%%%%%%%%
%%%%%%%%%%%%%%%%%%%%%%%%%%%%%%%%%%%%%%%%%%%%%%%%%%%%%%%%%%%%

%%% BEGIN PAPER

\setcounter{section}{-1}
\subsection*{Contents}

\begin{enumerate}
  \item[0.] Introduction 
  \item[1.] Generalities about discrepancy coefficients
  \item[2.] Minimal discrepancies and Shokurov's conjecture
  \item[3.] Proof of the main result (Theorem 1)
  \item[4.] Minimal discrepancies of log-terminal threefold singularities
\end{enumerate}

%%%%%%%%%%%%%%%%%%%%%%%%%%%%%%%%%%%%%%%%%%%%%%%%%%%%%%%%%%%%
%%%%%%%%%%%%%%%%%%%%%%%%%%%%%%%%%%%%%%%%%%%%%%%%%%%%%%%%%%%%

\section{Introduction}

  Let $Y$ be a normal, \Q-Gorenstein projective variety, and let $f:X \to Y$
  be a resolution of singularities. The discrepancy divisor $\Delta=K_X-f^*K_Y$
  plays a key role in the geometry of $Y$.  For example, the singularities
  allowed on a minimal (resp. on a canonical) model of $Y$ are defined in terms
  of $\Delta$.  Also, effective results for global generation of linear systems
  on singular threefolds (cf. \cite{elm}) involve certain coefficients of
  $\Delta$.

  There are many difficult conjectures, and several important results (at least
  in dimension $\leq 3$), regarding the discrepancy coefficients of $Y$ (i.e.,
  the coefficients of $\Delta$). In this paper we study a special case of
  the following problem:

\vspace{3pt}

  \noindent
  {\bf Shokurov's conjecture.} (cf. \cite{sho}, \cite{kollar})
     \emph{ If $\dim(Y) = n$, and $y \in Y$ is a singular point, then
     $\md_y(Y) \leq n-2$}. 

  The \emph{minimal discrepancy} of $Y$ at $y$, $\md_y(Y)$, is defined as
  $$ \md_y(Y) = \inf \{ \ord_F(\Delta) \mid f(F) = \{ y \} \ \ \}. $$

\vspace{3pt}

  The main result of this paper is an elementary computation of $\md_y(Y)$ for
  a large class of hypersurface singularities:

  \begin{Theorem}
    Assume that the germ $(Y,y)$ is analytically equivalent to a hypersurface
    singularity $(Y',0) \subset ({\mathbb{A}}_{\,\C}^{n+1},0)$, given by
      $$ Y'=\{ (y_1, \dots, y_{n+1}) \mid G(y_1, \dots, y_{n+1}) = 0 \} ;\quad
		    G(0, \dots , 0) = 0.	$$
    For an $n$-tuple $(a_1, \dots, a_n)$ of positive integers, write
    $G(t^{a_1}u_1, \dots, t^{a_n}u_n, t) =
  \\
    t^A \phi(u_1, \dots, u_n) +
      t^{A+1} \psi(u_1, \dots, u_n, t)$ with $\phi(u_1, \dots, u_n) \neq 0$.
    Note that $\phi$ is always a polynomial of degree at most $A$, even if
    $G$ is a power series. Assume that $\phi$ has at least one irreducible
    factor with exponent 1 in its factorization. 

    Then $\md_y(Y) \leq d$, where $d=(a_1+\cdots+a_n)-A$.    
  \end{Theorem}

  This criterion applies, for example, to hypersurface singularities of rank
  at least $2$ (if the singularity $(Y,0)$ is defined by $G(y_1,\dots,y_{n+1})
  =0$, then we define its rank as the rank of the quadratic part of $G$). It
  applies also to terminal (and, more generally, $cDV$) singularities in
  dimension 3. Shokurov's conjecture for terminal threefolds was proved
  by D. Markushevich in \cite{mark}, using the language of toric geometry,
  Newton diagrams, admissible weights, etc., and using the fact that the 
  singularities are isolated. Our proof shows that the result is completely
  elementary, and works for non-isolated singularities as well.

\vspace{5pt}

  The paper is organized as follows: In \S1 we discuss discrepancy coefficients
  in general. Everything in this section is well-known to the experts; we wrote
  it mainly to fix our notations and terminology.  We discuss in some detail
  the invariance of certain definitions under analytic equivalence of germs;
  we couldn't find a satisfactory reference in the literature.  (N. Mohan Kumar
  pointed out to us that the matter is not completely trivial.) In \S2 we
  introduce minimal discrepancies and prove some easy reductions of Shokurov's
  conjecture. In \S3 we prove Theorem 1, and in \S4 we carry out the
  computations for log-terminal threefold singularities.

\vspace{5pt}

  We would like to express our gratitude to L.~Ein, P.~Ionescu, R.~Lazarsfeld,
  K.~Matsuki, and N.~Mohan Kumar; our many conversations were very useful.

%%%%%%%%%%%%%%%%%%%%%%%%%%%%%%%%%%%%%%%%%%%%%%%%%%%%%%%%%%%%
%%%%%%%%%%%%%%%%%%%%%%%%%%%%%%%%%%%%%%%%%%%%%%%%%%%%%%%%%%%%

\section{Generalities about discrepancy coefficients}

  In this section we recall several definitions and results regarding
  discrepancy coefficients, cf. \cite{reid}, \cite{ckm}, \cite{kmm}.

\vspace{5pt}

{\bf (1.1)}
  Let $f:X \to Y$ be a birational morphism of normal projective varieties
  of dimension $n$ over \C. A prime Weil divisor $F \subset X$ is
  \emph{$f$-exceptional} if $\dim f(F) \leq n-2$. The closed subset
  $f(F) \subset Y$ is called the \emph{center} of $F$ on $Y$. More generally,
  a \Q-Weil divisor $D=\sum a_j F_j$ is $f$-exceptional if all the irreducible
  components $F_j$ are $f$-exceptional. Let $\Exc(f) = \union \{ F_j \mid
  \textup{$F_j \subset X$ a prime $f$-exceptional divisor} \}$; thus $D$ is
  $f$-exceptional if and only if $\Supp(D) \subset \Exc(f)$.

\vspace{5pt}

{\bf (1.2)}
  Choose a canonical divisor $K_Y$ on $Y$. Assume that $Y$ is \Q-Gorenstein,
  with global index $r$; i.e., $m K_Y$ is Cartier for some integer $m \geq 1$,
  and $r$ is the smallest such integer. Then we can define a \Q-divisor 
  $f^*K_Y$ on $X$ by $f^*K_Y = \tfrac{1}{r} f^*(r K_Y)$. On the other hand,
  there is a unique canonical divisor $K_X$ on $X$ such that the \Q-divisor 
  $\Delta = K_X - f^*K_Y$ is $f$-exceptional. ($K_X$ is obtained as follows:
  let $\omega$ be a rational differential $n$-form on $Y_{\textup{reg}}$, the
  smooth locus of $Y$; then $f^*\omega$ extends uniquely to a rational form on
  $X$, which we still denote by $f^*\omega$. If $\omega$ is chosen such that 
  $K_Y = \divisor_Y(\omega)$, then $K_X = \divisor_X(f^*\omega)$.) The divisor
  $\Delta = K_X - f^*K_Y$ is called the \emph{discrepancy divisor} of $f$. Note
  that $K_Y$ varies in a linear equivalence class on $Y$, and correspondingly
  $K_X$ varies in its own linear equivalence class on $X$; however, $\Delta$ is
  uniquely determined by $f$: indeed, if $\omega' = \phi \, \omega$ on
  $Y_{\textup{reg}}$, for some rational function $\phi \in \C(Y_{\textup{reg}})
  = \C(Y)$, then $f^*\omega' = (f^*\phi)(f^*\omega)$, and
  $\divisor_X(f^*\phi) = f^* \divisor_Y(\phi)$.

\vspace{5pt}

{\bf (1.3)}
  Write $\Delta = \sum a_j F_j$; the rational numbers $a_j$ are called
  \emph{discrepancy coefficients}. Now consider another birational morphism
  $f' : X' \to Y$ (with $X'$ a normal projective variety of dimension $n$).
  ${f'}^{-1} \circ f$ is a birational map $g: X \, \cdots \! \to X'$. Let
  $F_j \subset X$ be an $f$-exceptional divisor which intersects the regular
  locus $\Reg(g)$ of $g$, and assume that $g$ is an isomorphism at the generic
  point of $F_j$; i.e., $\overline{g(F_j \inters \Reg(g))}$ is a \emph{divisor}
  $F_j'$ on $X'$. Then $F_j'$ is an $f'$-exceptional divisor; in fact, $F_j$
  and $F_j'$ have the same center on $Y$, $f(F_j)=f'(F_j')$. Moreover, if
  $a_j'$ is the coefficient of $F_j'$ in $\Delta' = K_{X'} - {f'}^* K_Y$, then
  $a_j' = a_j$. This is seen by resolving the indeterminacies of the map $g$
  to a morphism $\tilde{g}:\tilde{X} \to X'$ ($\tilde{X}$ being obtained after
  a finite sequence of blowing-ups from $X$, cf. \cite[p.144, Consequence (1)
  of Corollary 1]{hironaka}) and calculating in two ways the discrepancy
  coefficient of $\tilde{F}_j$ = proper transform of $F_j$ on $\tilde{X}$.
  Thus, in fact, the discrepancy coefficient (and the center on $Y$) depends
  only on the discrete valuation of the rational function field $\C(Y)$
  determined by $F_j$ (or $F_j'$).

\vspace{5pt}

{\bf (1.4.)}
  Conversely, let $v$ be any discrete valuation of $\C(Y)$. Then $v$ is
  associated to a certain divisor $F^0 \subset X^0$ for some birational
  morphism $f^0 : X^0 \to Y$; in fact, by Hironaka's embedded resolution of
  singularities, if we start with \emph{any} birational morphism $f:X \to Y$
  as before, we can find a suitable $f^0$ with $X^0$ smooth, $\Exc(f^0)$ a
  divisor with normal crossings, and $X^0$ obtained from $X$ by a finite
  sequence of blowing-ups along smooth centers. Then $f^0(F^0) \subset Y$
  depends only on $v$; this closed subset is called the \emph{center of $v$ on
  $Y$}. $v$ is \emph{$Y$-exceptional} if this center has dimension at most
  $n-2$, and in this case $v$ has a well-defined discrepancy coefficient with
  respect to $Y$.

\vspace{5pt}

{\bf (1.5)}
  Let $f:X \to Y$ be as before, and let $F_j \subset X$ be $f$-exceptional.  The
  computation of the discrepancy coefficient $a_j$ is local on $X$; i.e., we
  may replace $Y$ with an open neighborhood of the generic point of $f(F_j)$,
  and $X$ with an open neighborhood of the generic point of $F_j$. From this
  point of view, the projectivity requirement is irrelevant. In particular, we
  may consider discrepancy coefficients for \emph{germs} $(Y,y)$ of algebraic
  varieties; one such coefficient is associated to each $Y$-exceptional
  discrete valuation of $\C(Y)$ whose center on $Y$ contains $y$.

  Moreover, the requirement that $X$ be normal is also irrelevant in some
  situations; for example, if $F_j$ is a Cartier divisor on $X$ (or at least on
  some open subset $U \subset X$ with $F_j \inters U \neq \emptyset$), then the
  generic point of $F_j$ has a \emph{nonsingular} open neighborhood in $X$, and
  we may replace $X$ with this neighborhood if we are interested only in the
  discrepancy coefficient of $F_j$.

\vspace{5pt}

{\bf (1.6)}
  {\bf Definition.} A projective variety $Y$ as before (i.e. normal,
    \Q-Gorenstein, $n$-dimensional) has \emph{only terminal (canonical,
    log-terminal, log-canonical) singularities} if all discrepancy coefficients
    of $Y$ are $>0$ (resp. $\geq 0$, $>-1$, $\geq -1$). Similarly, $Y$ is
    terminal (canonical, etc.) at a point $y$, or the germ $(Y,y)$ is terminal
    (etc.), if all discrepancy coefficients of discrete valuations with center
    containing $y$ are $>0$ (resp. $\geq 0$, etc.)

  \begin{Proposition}  \label{p2}
   \textup{(cf. \cite[Proposition 6.5]{ckm})}
    Let $f:X \to Y$ be a proper birational morphism, with $X$ smooth and
    $\Exc(f)$ with only normal crossings. Let $\Delta = K_X - f^*K_Y =
    \sum a_j F_j$, and let $\alpha = \min \{ a_j \}$.

    If $-1 \leq \alpha \leq 1$, then \emph{all} the discrepancy coefficients
    of $Y$ are $\geq \alpha$ (even for those discrete valuations of $\C(Y)$
    which are $Y$-exceptional but do not correspond to divisors on $X$).
  \end{Proposition}

  \vspace{3pt}

  In particular, to check whether $Y$ (or a germ $(Y,y)$) is terminal (etc.),
  it suffices to examine the discrepancy coefficients of a single resolution
  of singularities $f$ as above.

  \vspace{3pt}

  We reproduce the proof here for the reader's convenience (cf. \cite{ckm});
  the same computation will be used again in (1.7) and in (2.1).

 \begin{proof}
    As explained in (1.4), it suffices to consider a single blowing-up of $X$
    along a smooth center $Z \subset X$. Let $h:X' \to X$ be this blowing-up,
    $F_j' = h^{-1} F_j$ (proper transform), and $F'$ = the exceptional divisor
    of $h$. Let $r = \codim_X(Z) \geq 2$. Since $\Exc(f) = \union \, F_j$ has
    only normal crossings, $Z$ is contained in at most $r$ of the divisors
    $F_j$; say $Z \subset F_1, \ldots, F_s$, $s \leq r$. Let $f' = f \circ h :
    X' \to Y$, and $\Delta' = K_{X'} - {f'}^*K_Y$; then 
   \begin{multline*}
       \Delta' = K_{X'} - h^* f^* K_Y = K_{X'} - h^*(K_X-\Delta) = \\
       = K_{X'} - h^* K_X + h^*(\sum a_j F_j)
       = (r-1) F' + \sum a_j F_j' + (\sum_{j=1}^s a_j) F';
   \end{multline*}
    the discrepancy coefficient of $F'$ is therefore $a' = (r-1) +
    (\sum_1^s a_j)$. If $\alpha \leq 0$, we have $\sum_1^s a_j \geq
    s \alpha \geq r \alpha$ (because $s \leq r$ and $\alpha \leq 0$), and
    therefore $a' \geq (r-1) + r \alpha \geq \alpha$ (because $r > 1$ and
    $\alpha \geq -1$). If $\alpha > 0$, then we get $a' \geq r-1 \geq 1 \geq
    \alpha$.
 \end{proof}

 \emph{Remarks.} \ \emph{1.} The condition $\alpha \leq 1$ can always be
   achieved, as follows: let $f:X \to Y$ be a resolution of singularities,
   as in the statement of the proposition; choose a smooth subvariety
   $T \subset X$ of codimension $2$, such that $T \notsubset \Exc(f)$; and
   replace $f$ with $f \circ g$, where $g:\tilde{X} \to X$ is the blowing-up
   of $X$ along $T$. The computation used in the proof of the proposition
   shows that the exceptional divisor of $g$ has discrepancy coefficient $1$
   relative to $Y$.

\vspace{3pt}

   \emph{2.} If $\alpha < -1$, then the infimum of all discrepancy
   coefficients relative to $Y$ is $-\infty$; see \cite[Claim 6.3]{ckm}.
   (We prove a more precise statement in \S2, Lemma \ref{notlc}.)

   In general, the infimum of all discrepancy coefficients is called the
   \emph{(total) discrepancy} of $Y$, notation: $\discrep(Y)$. Thus
   $\discrep(Y) = -\infty$ if $Y$ is not log-canonical; if $Y$ \emph{is}
   log-canonical, then $-1 \leq \discrep(Y) \leq 1$, and $\discrep(Y)$ can be
   calculated by examining a single resolution of singularities $f:X \to Y$
   as in the proposition.

\vspace{3pt}

   \emph{3.} If $0 \leq \alpha \leq 1$, the proof shows that every 
   $Y$-exceptional discrete valuation of $\C(Y)$, other than those associated
   to the exceptional divisors of $f$, has discrepancy coefficient $\geq 1$.

\vspace{5pt}

{\bf (1.7)}
  Let $(Y,y)$ be an algebraic germ, as before, and let $(Y^{an},y)$ be the
  corresponding analytic germ; note that $Y$ normal and irreducible 
  $\implies$ $Y^{an}$ normal and irreducible. Also, $Y$ \Q-Gorenstein
  $\implies$ $Y^{an}$ \Q-Gorenstein. The theory of discrepancy divisors,
  discrepancy coefficients, terminal singularities, etc., can be developed
  in parallel in the category of germs of Moishezon analytic spaces; the
  results discussed so far are identical in the two categories.

  An interesting question arises when we try to compare the discrepancy
  coefficients for $(Y,y)$ and $(Y^{an},y)$. For example, is it true that
  $(Y,y)$ is terminal if and only if $(Y^{an},y)$ is terminal? (If this is
  true, then ``terminal'' depends only on the analytic equivalence class of an
  algebraic germ.) In general, the field of meromorphic functions of $Y^{an}$,
  $\Mcal(Y^{an})$, has many discrete valuations which vanish identically on
  $\C(Y)$; therefore the question is non-trivial.

   The answer is given by the following observation:

\begin{Proposition}		\label{dcan}
  Let $f:X \to (Y,y)$ be a proper birational morphism with $X$ smooth and
  $\Exc(f)$ with normal crossings.  Let $\{F_j\}_{j \in J}$ be the
  $f$-exceptional divisors on $X$, and let $\Delta = \sum a_j F_j$.

  Then the set of \emph{all} discrepancy coefficients of $(Y,y)$ is completely
  determined by the following combinatorial data:
    \begin{itemize}
      \item[(1)] The finite set $J$;
      \item[(2)] The rational numbers $a_j$ (one for each $j \in J$); and
      \item[(3)] For each subset $I \subset J$, the logical value of 
        $\displaystyle{\operatornamewithlimits{\bigcap}_{j \in I} F_j}
        \neq \emptyset$ \textsc{(true \textnormal{or }false)}.
    \end{itemize}
\end{Proposition}

  This observation (and its proof below) is valid in the algebraic as well as
  in the analytic case. In particular, the set of all ``algebraic'' and the set
  of all ``analytic'' discrepancy coefficients of $(Y,y)$ coincide.  (We may
  start with the same algebraic resolution $f:X \to (Y,y)$ in the analytic
  category, as $f^{an}:X^{an} \to (Y^{an}, y)$; then the initial combinatorial
  data for $f^{an}$ is the same as for $f$.)

\begin{proof}
  Let $v$ be a $Y$-exceptional discrete valuation of $\C(Y)$ with center
  containing~${y}$. By \cite[Main Theorem II]{hironaka}, there exists a finite
  succession of blowing-ups $f_i : X_{i+1} \to X_i$ along $Z_i \subset X_i$,
  where $0 \leq i < N$ and $X_0 = X$, with the following properties:
    \begin{itemize}
       \item[(i)] $v$ corresponds to a divisor on $X_N$;
       \item[(ii)] $Z_i$ is smooth and irreducible; and
       \item[(iii)] If $E_0=\Exc(f)$, and $E_{i+1}=f_i^{-1}(E_i)_{\text{red}}
       \union f_i^{-1}(Z_i)_{\text{red}}$, $0\leq i<N$, then $E_i$ has only
       normal crossings with $Z_i$.
    \end{itemize}
  (Recall what this means, from \cite[Definition 2]{hironaka}: at each point 
  $x \in Z_i$ there is a regular system of parameters of $\Ocal_{X_i,x}$, 
  say $(z_1, \dots, z_n)$, such that each component of $E_i$ which passes
  through $x$ has ideal in $\Ocal_{X_i,x}$ generated by one of the $z_j$, and
  the ideal of $Z_i$ in $\Ocal_{X_i,x}$ is generated by some of the $z_j$.)

  Let $f_1 : X_1 \to X$ be the blowing-up along a smooth irreducible subvariety
  $Z \subset X$, of codimension $r \geq 2$, such that $\Exc(f) =
  \union_{j \in J} F_j$ has only normal crossings with $Z$. Say $Z \subset F_j$ 
  if and only if $j \in \{j_1, \dots, j_s \}$; $s \leq r$.

  Considering $g = f \circ f_1 : X_1 \to Y$, we get a new element $j'$ added
  to $J$, $J_1 = J \union \{j'\}$, where $j'$ corresponds to the exceptional
  divisor $F'$ of $f_1$. The corresponding number is $a_{j'}=(r-1) + (a_{j_1}+
  \cdots + a_{j_s})$. Since the $F_j$ have only normal crossings with $Z$, the
  ``intersection data'' for $J_1$ is completely determined by the data for $J$,
  plus the following combinatorial data for $Z$:
    \begin{itemize}
       \item[(4)] For each $I \subset J$, the non-negative integer $d_I =
         \dim \left( Z \inters \left[
  \displaystyle{\operatornamewithlimits{\bigcap}_{j \in I} F_j}
                               \right] \right) $.
    \end{itemize}
  (Note that this collection of data contains, in particular, the codimension
  $r$ of $Z$, in the form $d_{\emptyset} = r$, and also the information about
  which $F_j$'s contain $Z$, in the form $Z \subset F_j \Leftrightarrow
  d_{\{j\}} = r$.)

  Finally, which such functions $\{d_I\}_{I \subset J}$ are possible is
  completely determined by the ``intersection data'' for $J$. Since every
  discrepancy coefficient of $Y$ is obtained after a finite number of such
  elementary operations on the combinatorial data (corresponding to a
  succession of blowing-ups along smooth centers), the result follows by
  induction.
\end{proof}

%%%%%%%%%%%%%%%%%%%%%%%%%%%%%%%%%%%%%%%%%%%%%%%%%%%%%%%%%%%%
%%%%%%%%%%%%%%%%%%%%%%%%%%%%%%%%%%%%%%%%%%%%%%%%%%%%%%%%%%%%

\section{Minimal discrepancies and Shokurov's conjecture}

{\bf (2.1)}
  {\bf Definition.} Let $(Y,y)$ be an algebraic or analytic germ (as always,
  we assume it is normal, \Q-Gorenstein, $n$-dimensional). The \emph{minimal
  discrepancy of $Y$ at $y$}, $\md_y(Y)$, is the infimum of all discrepancy
  coefficients of discrete valuations of $\C(Y)$, resp. $\Mcal(Y^{an})$,
  whose center on $Y$ is $y$.

  \begin{Lemma}  \label{notlc}
    If $(Y,y)$ is not log-canonical at $y$, then $\md_y(Y) = -\infty$.
  \end{Lemma}

  \begin{proof}
    Let $f \colon X \to (Y,y)$ be a resolution of singularities with $\Exc(f)$
    having only normal crossings. Let $F_j \subset X$ be an
    $f$-exceptional divisor with $y \in f(F_j)$ and discrepancy coefficient
    $a_j < -1$. Since $f^{-1}(y)$ is a union of $f$-exceptional divisors, and
    $F_j$ meets $f^{-1}(y)$, there is at least one exceptional divisor $F_i$
    with $f(F_i) = \{ y \} $ and $F_i \inters F_j \neq \emptyset$. We may
    assume that $F_i$ and $F_j$ are distinct (if $F_j \subset f^{-1}(y)$ and
    it is the only component of the fiber, we may blow up $X$ at a point of
    $F_j$; then take the exceptional divisor of this blowing-up in place of
    $F_i$, and the proper transform of $F_j$ in place of $F_j$). Set $Z = F_i
    \inters F_j$; then $Z$ is a smooth subvariety of codimension $2$ in $X$,
    and is not contained in any other exceptional divisor. Let $a_i$ be the
    discrepancy coefficient of $F_i$.

    Let $g: \tilde{X} \to X$ be the blowing-up of $X$ along $Z$. Let $F'$ be
    its exceptional divisor, with discrepancy coefficient $a'$ relative to $Y$.
    Then $a' = 1 + a_j + a_i$ (see the proof of Proposition \ref{p2} in \S1).
    Moreover, $F'$ has center $\{ y \}$ on $Y$, and intersects the proper
    transform $F_j'$ of $F_j$ on $\tilde{X}$ (which has discrepancy
    coefficient $a_j' = a_j$ relative to $Y$).

    Note that $a_j < 1 \implies a' < a_i$. In fact, since all the discrepancy
    coefficients of $(Y,y)$ are integer multiples of $\frac{1}{r}$ (if $r$ is
    the index of $K_Y$ at $y$), we see that $a' \leq a_i - \frac{1}{r}$.
    Therefore the proof may be completed by induction.
  \end{proof}

\vspace{5pt}

{\bf (2.2)}
  Recall the statement of Shokurov's conjecture from the Introduction. The
  lemma we have just proved shows that the conjecture is true for
  non-log-canonical singularities.

  Shokurov's conjecture is trivially true for curves (there are no singular
  normal points in dimension 1). It is also true in dimension 2: if $(Y,y)$
  is a normal singularity and $f:X \to (Y,y)$ is the \emph{minimal}
  desingularization, then \emph{all} the coefficients of $\Delta=K_X-f^*K_Y$
  are $\leq 0$ (see, for example, \cite[Lemma 1.4]{elm}).

\vspace{5pt}

{\bf (2.3)}
  The following lemma shows that the conjecture can be reduced to the case of
  singularities of index one:

  \begin{Lemma}   \label{indexone}
    Let $\varphi : Y' \to Y$ be a finite morphism of normal, \Q-Gorenstein
    varieties.  Assume that $\varphi$ is \'{e}tale in codimension one.
    Let $y'$ be a point of $Y'$, and $y = \varphi(y')$.

    Then $\md_y(Y) \leq \md_{y'}(Y')$.
  \end{Lemma}

  In particular, if $(Y,y)$ has index $r$, then there exists a $\varphi:Y'\to Y$
  as in the lemma, with $Y'$ having index one (the ``index-one cover'', cf.
  \cite[Definition 6.8]{ckm}).  Thus it would suffice to prove Shokurov's
  conjecture for singularities of index one.

  \begin{proof}   (cf. \cite[proof of Lemma 2.2]{elm})
    Let $f' : X' \to Y'$ be a resolution of singularities of $Y'$ such that 
    $\md_{y'}(Y')=\ord_{F'}(\Delta')\geq-1$, where $\Delta'=K_{X'}-{f'}^*K_{Y'}$
    and $F'$ is a divisor on $X'$ with $f'(F')=\{y'\}$.  (If $\md_{y'}(Y') =
    -\infty$, let $\alpha < -1$ be a rational number, and choose $f', F'$ such
    that $f'(F') = \{y'\}$ and $\ord_{F'}(\Delta') < \alpha$.)

    Let $f:X \to Y$ be a resolution of singularities of $Y$. By blowing up
    $X$, then $X'$, if necessary, we may assume that $\psi = f^{-1} \circ
    \varphi \circ f' : X' \to X$ is a morphism and that $\psi(F')$ is a
    divisor $F \subset X$. Let $\Delta = K_X - f^*K_Y$ and $a=\ord_F(\Delta)$.

    Let $t$ be the ramification index of $\psi$ along $F'$. Then:
    \begin{align*}
      K_{X'} &= \psi^*K_X + (t-1)F' + \textup{ other terms }		\\
             &= \psi^*(f^*K_Y + aF + \textup{ other terms }) + (t-1)F' +
             			\textup{ other terms }			\\
             &= \psi^*f^*K_Y + (ta+t-1)F' + \textup{ other terms }	\\
             &= {f'}^*\varphi^*K_Y  + (ta+t-1)F' + \textup{ other terms }  \\
             &= {f'}^*K_{Y'}  + (ta+t-1)F' + \textup{ other terms }	
    \end{align*}
    (note that $\varphi^*K_Y = K_{Y'}$, because $\varphi$ is \'{e}tale in
    codimension one).  Therefore $\ord_{F'}(\Delta') = ta+t-1$.

    If $\ord_{F'}(\Delta') \geq -1$, we get
    $\ord_F(\Delta) = a \leq \ord_{F'}(\Delta') = \md_{y'}(Y')$; indeed,
    $t \geq 1$, and therefore $ta \leq ta + (t-1)(1+\ord_{F'}(\Delta')) =
    t\ord_{F'}(\Delta')$.
    
    If $\ord_{F'}(\Delta') < \alpha < -1$, then $\ord_F(\Delta) = a \leq
    \tfrac{1}{t} \ord_{F'}(\Delta') < \tfrac{1}{t}\alpha$, with $1 \leq t
    \leq \deg(\varphi)$ and $\alpha$ an arbitrarily negative rational number.

    Since $f(F) = f\psi(F') = \varphi f'(F') = \varphi(\{y'\}) = \{y\}$, the
    lemma is proved.
  \end{proof}

\vspace{5pt}

{\bf (2.4)}
  Finally, we show that $\md_y(Y)$ is an analytic invariant.  In fact, we show
  that the set of all discrepancy coefficients for divisors with center
  $\{y\}$ on $Y$ is the same in the algebraic and in the analytic category:

  \begin{Proposition}		\label{mdan}
    Let $f:X \to (Y,y)$ be a resolution of singularities, as in Proposition
    \ref{dcan}. Then the set of all discrepancy coefficients for divisors with
    center $\{y\}$ on $Y$ is completely determined by the combinatorial data
    (1), (2), (3) in Proposition~\ref{dcan}, plus:
        \begin{itemize}
            \item[(3+)] For each $j \in J$, the logical value of
            ``$f(F_j)=\{y\}$'' \textsc{(true \textnormal{or} false)}.
        \end{itemize}
  \end{Proposition}

  \begin{proof}
      Let $f_1:X_1 \to X$ be the blowing-up of a smooth subvariety $Z\subset X$,
      as in the proof of Proposition \ref{dcan}. Put $g=f \circ f_1$, and let
      $F'$ be the exceptional divisor of $f_1$. Then $[g(F')=\{y\}]
      \Leftrightarrow [f(F_j)=\{y\}$ for at least one of the $F_j$'s containing
      $Z$]. Indeed, if $Z \subset F_j$ and $f(F_j)=\{y\}$, then $g(F') = f(Z)
      \subset \{y\}$, so that in fact $g(F')=\{y\}$. Conversely, $g(F') =
      \{y\} \implies Z \subset f^{-1}(y)$. As $Z$ is irreducible and
      $f^{-1}(y)$ is a union of divisors $F_j$ with $f(F_j)=\{y\}$, $Z$ must
      be contained in at least one such $F_j$.

      Therefore the ``extended'' combinatorial data for $g$ (including the
      information in (3+)) can be obtained from the ``extended'' combinatorial
      data for $f$.  The conclusion follows by induction.
  \end{proof}
  
%%%%%%%%%%%%%%%%%%%%%%%%%%%%%%%%%%%%%%%%%%%%%%%%%%%%%%%%%%%%
%%%%%%%%%%%%%%%%%%%%%%%%%%%%%%%%%%%%%%%%%%%%%%%%%%%%%%%%%%%%

\section{Proof of the main result (Theorem 1)}

{\bf (3.1)}
  Recall the statement of Theorem 1 from the Introduction. By (2.4), we may
  assume that $Y$ is the hypersurface $G=0$ in $\A^{n+1}$, with $y=0$.
  For convenience, denote $\A^{n+1}$ by $V$; thus $Y \subset V$.
  Let $U = \A^{n+1}$; write the coordinates in $V$ as $(y_1, \dots, y_{n+1})$,
  and the coordinates in $U$ as $(u_1, \dots, u_n, t)$. 

Let $f:U \to V$ be the birational morphism defined by $y_{n+1}=t; 
y_i = t^{a_i} u_i, i=1, \dots, n$. Let $E \subset U$ be the hyperplane
$(t=0)$; then $\Exc(f)=E$.

\vspace{5pt}

{\bf (3.2)}
  Let $\bar{Y} \subset U$ be the proper transform of $Y$ by $f$, $\bar{f} :
  \bar{Y} \to Y$ the restriction of $f$ to $\bar{Y}$, and $\bar{E} =
  E |_{\bar{Y}}$ (as a Cartier divisor). By hypothesis, $f^*Y = \bar{Y} + AE$,
  and $\bar{E}$ has equation $\phi(u_1, \dots, u_n) = 0$ in $E \cong \A^n$.
  Since $\phi$ has at least one irreducible factor with exponent 1, $\bar{E}$
  has at least one irreducible component with multiplicity one: $\bar{E}=F_1+
  \cdots $. As explained in (1.5), since $\bar{Y}$ is smooth in a neighborhood
  of the generic point of $F_1$, and $F_1$ is the exceptional divisor which will
  produce the desired discrepancy coefficient, we need not worry about the
  normality of $\bar{Y}$. 

\vspace{5pt}

{\bf (3.3)}
  Take $\omega = dy_1 \wedge \cdots \wedge dy_{n+1}$ on $V$; then 
  $f^*\omega = t^{a_1+\cdots+a_n}du_1 \wedge \cdots \wedge du_n \wedge dt$ on
  $U$, so that $K_U - f^*K_V = (a_1 + \cdots a_n)E$.

  The adjunction formula gives $K_Y = K_V + Y |_Y$ and $K_{\bar{Y}} = K_U + 
  \bar{Y} |_{\bar{Y}}$. Therefore we have:
    \begin{align*}
       K_{\bar{Y}} - \bar{f}^* K_Y &= (K_U+\bar{Y}) |_{\bar{Y}} -
     			  \bar{f}^*(K_V+Y|_Y)				\\
     	  &= (K_U - f^*K_V + \bar{Y} - f^*Y) |_{\bar{Y}}		\\
     	  &= ((a_1 + \cdots + a_n) E - AE ) |_{\bar{Y}}			\\
     	  &= d \bar{E} = d F_1 + \cdots .
  \end{align*}
  (Recall that $d = (a_1 + \cdots + a_n) - A$.)

  Thus the discrepancy coefficient of $F_1\subset\bar{Y}$ with respect to $Y$
  is equal to $d$. Since $\bar{f}(F_1)=\{y\}$, Theorem 1 is proved.    \qed

\vspace{5pt}
{\bf (3.4) Example.}
  Let $(Y,y)$ be a singular germ of multiplicity 2. That is, we may assume that
  $Y$ is a hypersurface in $\A^{n+1}$ given by an equation $G=0$, with $y=0$,
  $G$ and all its first-order partial derivatives at 0 equal to zero, and some
  second-order partial derivative of $G$ at 0 not equal to zero.

  If $(Y,0)$ has rank at least 2, then $\md_y(Y) \leq n-2$ (as predicted by
  Shokurov's conjecture).  Indeed, consider the usual blowing-up of $V$ at 0;
  that is, take $a_1=\cdots=a_n=1$.  The hypothesis means that $A=2$, and ---
  after a linear change of parameters, if necessary --- $\phi(u_1,\dots,u_n) =
  u_1^2 + \cdots + u_r^2$, where $r \geq 2$ is the rank of the singularity.
  Thus $d=(a_1+\cdots+a_n)-A=n-2$ in this case, and $\phi$ is irreducible (if
  $r \geq 3$), resp.  a product of two distinct irreducible (linear) factors,
  if $r = 2$.

%%%%%%%%%%%%%%%%%%%%%%%%%%%%%%%%%%%%%%%%%%%%%%%%%%%%%%%%%%%%
%%%%%%%%%%%%%%%%%%%%%%%%%%%%%%%%%%%%%%%%%%%%%%%%%%%%%%%%%%%%

\section{Minimal discrepancies of log-terminal threefold singularities}

  Let $(Y,y)$ be a three-dimensional log-terminal singularity. In this section
  we will show that $\md_y(Y) \leq 1$ (so that Shokurov's conjecture is true
  in this case).

\vspace{5pt}

{\bf (4.1)}
  As shown in (2.3), we may assume that $(Y,y)$ has index one.  Then $(Y,y)$ is
  canonical (the index-one cover of a log-terminal singularity is again
  log-terminal, by Propositon \ref{p2}, and therefore canonical).

  In this case, M.~Reid \cite[Theorem 2.2]{reid} proved that either $(Y,y)$
  is a $cDV$ point (see below), or there exists a proper birational morphism
  $f:Y'\to Y$ with $f^* K_Y=K_{Y'}$ and $f^{-1}(y)$ containing at least one
  prime divisor of $Y'$.  Of course, in the latter case we have $\md_y(Y) = 0$.
  There only remains to consider the case when $(Y,y)$ is a compound Du Val
  ($cDV$) point; that is, $(Y,y)$ is analytically equivalent to a hypersurface
  singularity at the origin $0 \in \A^4$, with equation $G=0$,
      $$ G(y_1, y_2, y_3, t) = f(y_1,y_2,y_3) + t g(y_1,y_2,y_3,t), $$
  where $f(y_1,y_2,y_3)=0$ defines a Du Val singularity (rational double point)
  of a surface at $0 \in \A^3$.

  To simplify notation, we write $\bby$ for $y_1,y_2,y_3$ and $\bbu$ for 
  $u_1,u_2,u_3$.

  By (2.4), we may assume that $(Y,y)$ \emph{is} the hypersurface $(G=0)\subset
  \A^4$, with $y=0$.  By Theorem 1, it suffices to find $a_1,a_2,a_3 \geq 1$
  such that 
$$ G(t^{a_1}u_1,t^{a_2}u_2,t^{a_3}u_3,t)=t^A \phi(\bbu)+t^{A+1}\psi(\bbu,t) $$
  with $\phi(\bbu) \neq 0$, $(a_1+a_2+a_3)-A = 1$, and $\phi$ having at least
  one irreducible factor with exponent one in its prime decomposition.

\vspace{5pt}

{\bf (4.2)}
  We will do a case-by-case analysis, according to the type of singularity;
  $(Y,0)$ is of type $cA_n$, $cD_n$, or $cE_n$, if the surface singularity
  $f(\bby)=0 \subset \A^3$ is of type $A_n$, $D_n$, or $E_n$. 

  In each case, $f(\bby)$ is completely known.  $g(\bby,t)$, on the other hand,
  is not.  Of course, we have $g(0,0,0,0)=0$, or else $(Y,0)$ would be a smooth
  point.  We will not make any other assumptions about $g$.

  Write $g = g_1 + g_2 + \cdots$, where $g_i$ is a homogeneous form of degree
  $i$, and $o_k(g) = g_k + g_{k+1} + \cdots \,\, (k \geq 1)$.

  A note on terminology: we distinguish between \emph{form} and
  \emph{polynomial}; for instance, a quadratic polynomial is the sum of a
  quadratic form, a linear form, and a constant term.  We say that a polynomial
  (or a form) \emph{contains} a certain monomial if the coefficient of the
  monomial in that polynomial is non-zero.  We say that a monomial
  \emph{contains} $y_1$ if that monomial is divisible by $y_1$.

\vspace{5pt}

{\bf (4.3) Case} $\mathbf{cA_n}$: \ \ $f(\bby) = y_1^2 + y_2^2 + y_3^{n+1} 
			\quad (n \geq 1)$.

  Then $(Y,0)$ is a singularity of multiplicity 2 and rank at least 2.  This
  case is therefore covered by the Example discussed in (3.4).

\vspace{5pt}

{\bf (4.4) Case} $\mathbf{cD_n}$: \ \ $f(\bby) = y_1^2 + y_2 y_3^2 + y_3^{n-1} 
			\quad (n \geq 4)$.

  If $g_1(\bby,t) \neq 0$, then the quadratic part of $G = f + t g$ is 
  $y_1^2 + t g_1(\bby,t)$.  If this quadratic part has rank at least 2, then
  the conclusion follows from (3.4).  If it has rank 1, i.e. if
  $y_1^2+tg_1(\bby,t)$ is the square of a linear form, then a linear change of
  variable, $y_1' = y_1 + \alpha y_2 + \beta y_3 + \gamma t$, transforms the
  equation $G=0$ into a similar one with $g_1(\bby,t)=0$.

  So we need to consider only the case $g_1=0$.  Note that a similar argument
  applies to singularities of type $cE_n$.

  Assume that $g_1=0$. Then put $a_1=2, a_2=a_3=1$; that is, put $y_1=t^2 u_1,
  y_2 = t u_2, y_3 = t u_3$.  We have:
      $$ G(t^2u_1,tu_2,tu_3,t) = t^3\phi(\bbu) + t^4\psi(\bbu,t),	$$
  where $\phi(\bbu) = u_2 u_3^2 + \delta_{n,4} u_3^3 + [\text{terms of degree
  $\leq 2$ in the $u_j$}]$; $\delta_{n,4}=1$ if $n=4$, otherwise $\delta_{n,4}
  = 0$.  (The terms of lower degree come from $t g_2(t^2u_1,tu_2,tu_3,t)$, with 
  $g_2$ --- the quadratic component of $g(\bby,t)$.  Note that not \emph{all}
  the terms in $tg_2$ contribute to $\phi(\bbu)$: as $y_1=t^2u_1$, the
  monomials in $g_2(\bby,t)$ which contain $y_1$ give rise to monomials
  containing $t$ to the fourth or higher power.)

  The proof in this case is complete, for $(a_1+a_2+a_3)-A = (2+1+1)-3 = 1$,
  and $\phi$ has at least one irreducible factor with exponent one (otherwise
  $\phi$ would have to be the cube of a linear polynomial in $\bbu$; that
  linear polynomial would have to contain $u_2$, because $\phi$ contains
  $u_2 u_3^2$, and then $\phi$, being the cube of that linear polynomial,
  would contain $u_2^3$, which is not the case).

\vspace{5pt}

{\bf (4.5) Case} $\mathbf{cE_6}$: \ \ $f(\bby) = y_1^2 + y_2^3 + y_3^4$.

  As in (4.4), we may assume that the linear part $g_1(\bby,t)$ of $g(\bby,t)$
  is equal to zero.

  In $g_2$ (the quadratic part of $g$), separate the monomials which contain
  $y_1$ from those that don't: $g_2(\bby,t) = y_1 L(\bby,t) + Q(y_2,y_3,t)$,
  where $L$ is a linear form and $Q$ is a quadratic form.

  Put $y_1 = t^2 u_1, y_2 = t u_2, y_3 = t u_3$; then
      $$ G(t^2u_1, tu_2, tu_3, t) = t^3\phi(\bbu) + t^4\psi(\bbu,t),	$$
  with $\phi(\bbu) = u_2^3 + Q(u_2,u_3,1)$.

  If $\phi$ is not the cube of a linear polynomial, then we complete the proof
  just as in (4.4).  However, in this case it might be that $\phi$ \emph{is} a
  perfect cube.  If this is so, then $y_2^3 + t Q(y_2,y_3,t)$ is the cube of a
  linear form in $y_2,y_3,t$.  A linear change of variable, 
  $y_2' = y_2 + \alpha y_3 + \beta t$, reduces the proof to the case $Q=0$.
  This argument stands also in the cases $cE_7$ and $cE_8$, discussed below.

  There only remains to consider the case $g_2(\bby,t)=y_1L(\bby,t)$, where
  $L$ is a linear form (possibly zero).  In this case put $a_1=a_2=2, a_3=1$,
  i.e. $y_1=t^2u_1, y_2=t^2u_2, y_3=tu_3$.  Then:
    \begin{gather*}
        G(\bby,t) = y_1^2 + y_2^3 + y_3^4 + t[y_1L(\bby,t) + o_3(g)], \quad
    					    \text{and}			\\
        G(t^2u_1,t^2u_2,tu_3,t) = t^4 \phi(\bbu) + t^5 \psi(\bbu,t),
    \end{gather*}
 where $\phi(\bbu)=u_1^2+u_3^4+[\text{terms of degree $\leq 3$ in the $u_j$}]$,
 and the expression in brackets does not contain $u_1^2$ (so that $u_1^2$
 doesn't cancel out from $\phi$).

  Note that $\deg(\phi)=4$, and $\phi$ cannot be the square of a quadratic
  polynomial in the $u_j$ (if it were, then $\phi$ would contain the mixed
  product $u_1 u_3^2$, because it contains $u_1^2$ and $u_3^4$ but no $u_1^4$;
  the monomial $u_1u_3^2$ could only arise from a monomial $t(y_1y_3^2 t^k)$ of
  $tg(\bby,t)$, $k \geq 0$; but $t(t^2u_1)(tu_3)^2t^k = t^{5+k}u_1u_3^2$, so
  $u_1u_3^2$ cannot be a monomial of $\phi$).  Therefore $\phi$ has an
  irreducible factor with exponent one, and the conclusion follows --- note
  that $(a_1+a_2+a_3)-A = (2+2+1)-4 = 1$.

\vspace{5pt}

{\bf (4.6) Case} $\mathbf{cE_7}$: \ \ $f(\bby) = y_1^2 + y_2^3 + y_2 y_3^3$.

  We may again assume that $g_1=0$, as in (4.4), and that $g_2=y_1L(\bby,t)$
  with $L$ a linear form (possibly zero), as in (4.5).

  Write $L(\bby,t) = L_1(y_1,y_2)+L_2(y_3,t)$, and $g_3(\bby,t) = C_1(\bby,t) +
  C_2(y_3,t)$, where $C_1$ and $C_2$ are cubic forms such that every monomial
  of $C_1$ contains $y_1$ or $y_2$.

  Put $a_1=a_2=2, a_3=1$; then
   \begin{gather*}
    \begin{split}
     G(\bby,t)=y_1^2+y_2^3+y_2 y_3^3+t [ y_1 L_1(y_1,y_2)&+ y_1 L_2(y_3,t) \\
                 &+ C_1(\bby,t)+C_2(y_3,t)+o_4(g)],\quad\text{and}
    \end{split}							  \\
      G(t^2u_1,t^2u_2,tu_3,t) = t^4 \phi(\bbu) + t^5 \psi(\bbu,t), \qquad\qquad
       \qquad\qquad
   \end{gather*}
  where $\phi(\bbu) = u_1^2 + u_1 L_2(u_3,1) + C_2(u_3,1)$.

  Note that $(a_1+a_2+a_3)-A = (2+2+1)-4 = 1$.

  If $\phi$ has degree 3 (i.e. if $C_2(y_3,t)$ contains $y_3^3$), then $\phi$
  has an irreducible factor with exponent one, because $\phi$ cannot be a
  perfect cube (it contains $u_1^2$ but no $u_1^3$); in this case the proof is
  complete.

  Otherwise $\phi$ has degree 2 (it contains $u_1^2$).  Then either $\phi$ has
  an irreducible factor with exponent one (and then the proof is complete), or
  else $\phi$ is the square of a linear polynomial.  In the latter case, 
  $y_1^2 + t y_1 L_2(y_3,t) + t C_2 (y_3,t)$ is a perfect square.  The
  (non-linear) change of variable $y_1' = y_1 + \tfrac{1}{2} t L_2(y_3,t)$
  transforms the equation $G=0$ into a similar one with $L_2=C_2=0$
  (the verification is straightforward).  Therefore we may assume that $G$ has
  the form:
  $$ G(\bby,t) = y_1^2 + y_2^3 + y_2 y_3^3 + t [ y_1 L_1(y_1,y_2) + 
                C_1(\bby,t) + o_4(g)],					$$
  where $L_1$ is a linear form, and $C_1$ is a cubic form such that every
  monomial of $C_1$ contains $y_1$ or $y_2$.  (The same argument carries over
  unchanged to the last case, $cE_8$.)

  Now put $a_1=3,a_2=2,a_3=1$; that is, $y_1=t^3u_1, y_2=t^2u_2, y_3=tu_3$.
  We have:
      $$ G(t^3u_1, t^2u_2, tu_3,t) = t^5 \phi(\bbu) + t^6 \psi(\bbu,t),	$$
  where $\phi(\bbu) = u_2 u_3^3 + u_2 p(u_3) + q(u_3)$; $u_2 p(u_3)$
  corresponds to the monomials of $C_1(\bby,t)$ of the form $y_2 y_3^k t^{2-k},
  k=0,1,2$ (all other monomials of $C_1$ produce at least a sixth power of $t$;
  recall that all monomials of the cubic form $C_1$ contain $y_1$ or $y_2$),
  and $q(u_3)$ corresponds to the monomials of $g_4$ of the form
  $y_3^k t^{4-k}, k=0, \dots, 4$.  Note that $\phi$ has degree exactly one as a
  polynomial in $u_2$, and therefore $\phi$ cannot be the square of another
  polynomial in the $u_j$.  Since $\phi$ has (total) degree 4, it must have an
  irreducible factor with exponent one.  As $(a_1+a_2+a_3)-A = (3+2+1)-5 = 1$,
  the proof is complete in this case.

\vspace{5pt}

{\bf (4.7) Case} $\mathbf{cE_8}$: \ \ $f(\bby) = y_1^2 + y_2^3 + y_3^5$.

  As in the previous case, we may assume that
      $$ G(\bby,t) = y_1^2 + y_2^3 + y_3^5 + t [ y_1 L_1(y_1,y_2) + 
                C_1(\bby,t) + o_4(g)],					$$
  where $L_1$ is a linear form and $C_1$ is a cubic form such that every
  monomial of $C_1$ contains $y_1$ or $y_2$.

  If we take $a_1=3, a_2=2, a_3=1$, i.e. $y_1=t^3 u_1, y_2=t^2 u_2, y_3=tu_3$,
  we get
      $$ G(t^3u_1, t^2u_2, tu_3, t) = t^5 \phi(\bbu) + t^6 \psi(\bbu,t)	$$
  with $\phi(\bbu) = u_3^5 + u_2 p(u_3) + q(u_3)$, where $p(u_3)$ and $q(u_3)$
  are exactly as in the previous case.

  $(a_1+a_2+a_3)-A = (3+2+1) - 5 = 1$; so the proof is complete if $\phi$ has
  an irreducible factor with exponent one.

  Since $\deg(\phi) = 5$, if $\phi$ does not have an irreducible factor with
  exponent one then $\phi = M^2 N^3$ for two linear polynomials $M=M(u_2,u_3)$
  and $N = N(u_2,u_3)$ (possibly equal).  In particular, if this is the case
  then $\phi$ cannot have degree exactly one as a polynomial in $u_2$, and
  therefore $p(u_3)=0$; i.e., $C_1(\bby, t)$ contains no monomials of the form
  $y_2 y_3^k t^{2-k}\,(k=0,1,2)$.  As every monomial of $C_1$ contains $y_1$ or
  $y_2$, this means that every monomial of $C_1$ actually contains $y_1$ or
  $y_2^2$.

  Now $\phi(\bbu) = u_3^5 + q(u_3) \text{ ($q$ of degree at most four) } =
  M^2(u_3) N^3(u_3)$.  If we write $g_4(\bby,t) = F_1(\bby,t) + F_2(y_3,t)$,
  with $F_1,F_2$ forms of degree 4 such that every monomial of $F_1$ contains
  $y_1$ or $y_2$, then $q(u_3)=F_2(u_3,1)$.  $u_3^5+q(u_3) = M^2(u_3)N^3(u_3)$
  means that $y_3^5+F_2(y_3,t) = \tilde{M}^2(y_3,t)\tilde{N}^3(y_3,t)$, with
  $\tilde{M}, \tilde{N}$ linear forms in $y_3,t$ ($M(u_3)=\tilde{M}(u_3,1)$,
  etc.)  A linear change of variable $y_3'=y_3 + \alpha t$ reduces $G$ to the
  form
      $$ G(\bby,t) = y_1^2 + y_2^3 + y_3^3 ( y_3 + at)^2 + t[ y_1 L(y_1,y_2) +
  	    C_1(\bby,t) + F_1(\bby,t) + o_5(g)],			    $$
  where $a \in \C$ (possibly $a=0$), every monomial of the cubic form $C_1$
  contains $y_1$ or $y_2^2$, and every monomial of the quartic form $F_1$
  contains $y_1$ or $y_2$.

  Put $a_1=3, a_2=a_3=2$; that is, $y_1=t^3 u_1, y_2=t^2u_2, y_3=t^2 u_3$. Then
    $$ G(t^3u_1, t^2u_2, t^2u_3,t) = t^6 \phi(\bbu) + t^7 \psi(\bbu,t)	$$
  with $\phi(\bbu) = u_1^2 + u_2^3 + [\text{terms of degree at most 2 in the
  $u_j$}]$, and $u_1^2$ is not among the terms inside the brackets.  Therefore
  $\phi$ has an irreducible factor with exponent one; as $(a_1+a_2+a_3)-A =
  (3+2+2) - 6 = 1$, the proof is complete in all cases.			\qed

\vspace{5pt}

{\bf (4.8)}
  \emph{Remarks.} \ \emph{1.} If $(Y,y)$ is a terminal threefold singularity of
  index one, then $\md_y(Y) = 1$ ($\md_y(Y) \geq 1$ because all the discrepancy
  coefficients of $Y$ at $y$ are positive integers).  On the other hand,
  Kawamata \cite{kaw} proved that the minimal discrepancy of a terminal
  threefold singularity of index $r \geq 2$ is $\frac{1}{r}$.

\vspace{3pt}

  \emph{2.} Our computations in \S4 seem related to those in \cite{mark}, except
  for the fact that Markushevich uses the toric language.  At first glance, it
  looks like his proof uses the assumption that the singularity is isolated;
  however, this is needed only to reduce the equation $G=0$ to various standard
  forms, and -- as our elementary computations in (4.3) -- (4.7) show --- this
  can be done without assuming the singularity is isolated.  This allows us to
  conclude that Shokurov's conjecture is true for canonical (and log-terminal)
  threefold singularities, rather than just for terminal singularities.

%%%%%%%%%%%%%%%%%%%%%%%%%%%%%%%%%%%%%%%%%%%%%%%%%%%%%%%%%%%%
%%%%%%%%%%%%%%%%%%%%%%%%%%%%%%%%%%%%%%%%%%%%%%%%%%%%%%%%%%%%

    %%% BEGIN REFERENCES 

    %%% END REFERENCES

%%% END PAPER

\end{document}